\documentclass[12pt]{article}  
\usepackage[final]{graphicx}
\usepackage[dvips]{color}
\usepackage{graphicx}

\def\Bbb{\bf} 
\newcommand\C{{ \Bbb C}}   %
\newcommand\D{{ \Bbb D}}   %
\newcommand\R{{\Bbb R}}    %
\newcommand\Z{{\Bbb Z}}    %
\newcommand\N{{\Bbb N}}    %
\newcommand\PP{{\Bbb P}}
\newcommand\HH{{\Bbb H}}
\newcommand\ds{{\displaystyle}}

  \newtheorem{theorem}{Theorem}[section]
  
  \newtheorem{lemma}[theorem]{Lemma}
  
  \newtheorem{corollary}[theorem]{Corollary}
  \newtheorem{remark}[theorem]{Remark}{\rmfamily}

\def\comment#1{ }

\newcommand{\dfrac}[2]{%
  \frac{\displaystyle{#1}}{\displaystyle{#2}}}

\begin{document}
\title{Affine Schwarz map for the hypergeometric differential equation}
\author{Ryoichi Kobayashi, Tatsuya Nishizaka, Shoji Shinzato and Masaaki Yoshida
}

\maketitle

\begin{abstract}
We propose an affine version of the Schwarz map for the hypergeometric differential equation, and study its image when the monodromy group is finite.

\end{abstract}


\tableofcontents

\section{Introduction} Consider the {\it hypergeometric differential equation} 
$$E(a,b,c): x(1-x)w''+\{c-(a+b+1)x\}w'-abw=0,$$
and define its {\it affine Schwarz map} by
$$asch:X=\C-\{0,1\}\ni x\longmapsto (u,v)=(u(x),v(x))\in W\ (\cong\C^2),$$
where $u$ and $v$ are linearly independent solutions of $E(a,b,c)$. 
Since the {\it Schwarz map} is defined by
$$sch:X=\C-\{0,1\}\ni x\longmapsto z=u(x):v(x)\in Z\ (\cong\PP^1=W\!-\!\{(0,0)\}/\C^\times),$$
where $\PP^1$ is the complex projective line, our map $asch$ can be considered to be an affine version of $sch$. Note that the monodromy of $sch$ is the projectivization of that of $asch$ (i.e. of the equation $E(a,b,c)$). 

\par\medskip\noindent
It seems that no one has ever studied the affine Schwarz map before.
\par\smallskip\noindent 
In this paper, we first give an expression of the image curve under $asch$ when the monodromy group is finite and the inverse of the Schwarz map is single-valued; some other cases are also studied. A particular case, when the projective monodromy group is dihedral, yields an identity that the product of two hypergeometric series equals 1. We next generalize this and obtain identities, with two integer parameters, that the product of two hypergeometric series are polynomials. Finally we discuss briefly the image curve when the projective monodromy group is a Fuchsian group.
\section{When the inverse of the Schwarz map is single-valued}
When the monodromy group of the equation $E(a,b,c)$ is finite, the image under $asch$ should be an irreducible algebraic curve $C$ in $W$. We want to give a defining equation of $C$, and an expression of the inverse map of $asch:X\to C\subset W$.

Let $G\subset GL(2,\C)$ be the monodromy group with respect to the pair $(u,v)$; the projectivization $\bar G\subset PGL(2,\C)$ of $G$ is the monodromy group with respect to the ratio $u:v$. When the group $G$ is of finite order, $\bar G$ is isomorphic to one of the polyhedral groups. 
If a polynomial $P$, homogeneous in two variables $(u,v)$, is a semi-invariant of $G$, i.e.
$$P((u,v)g)=C_gP(u,v),\quad g\in G, \ C_g\in\C,$$
then the function $f(x)=P(u(x),v(x))$ on $X$ is the product of a rational function and fractional powers of $x$ and $1-x$:
$$ f(x)=(\mbox{a rational function in $x$})\cdot x^*\cdot(1-x)^*,$$
since any 1-dimensional representation of the fundamental group $\pi_1(X)$ of $X$ can be realized by the product of powers of $x$ and $1-x$, and since the singularities of the equation $E$ are regular.

When the inverse of the Schwarz map is single-valued, we give three identities of the form
$P(u(x),v(x))=f(x),$ from which the equation of the image $C$ and the inverse of $asch$ can be readily derived. Moreover, from each group of Schwarz's list, we choose a simple one, and do the same.
\par\medskip\noindent
Throughout this paper, we make use of the Kummer solutions 
$$u=u(x)=F(a,b,c;x),\quad v=v(x)=x^{1-c}F(a-c+1,b-c+1,2-c;x)$$
around the origin, where
$$F(a,b,c;x)=\sum_{n=0}^\infty\frac{(a,n)(b,n)}{(c,n)(1,n)}x^n,\quad 
(a,n)=a(a+1)\cdots(a+n-1).$$

\subsection{Preliminaries and notation}
If the parameters $(a,b,c)$ are real, the image of the upper part $X_+=\{x\in X\mid \Im x\ge0\}$ of $X$ under $sch$ of the Schwarz map of the equation $E(a,b,c)$ is a triangle with angles
$$|\mu_0|\pi,\ |\mu_1|\pi,\ |\mu_\infty|\pi, \quad\mbox{where}\quad \mu_0=1-c,\ \mu_1=c-a-b,\ \mu_\infty=a-b,$$
at the vertices $sch(0), sch(1),$ and $sch(\infty)$, respectively. Thus the $sch$ has a single-valued inverse only when
$$ |\mu_0|=\frac1{k_0},\quad |\mu_1|=\frac1{k_1},\quad |\mu_\infty|=\frac1{k_\infty},\quad k_0,k_1,k_\infty\in\{2,3,\dots,\infty\}.$$
If the monodromy group is finite, then $1/k_0+1/k_1+1/k_\infty>1$ must hold; this implies 
$$\{k_0,k_1,k_\infty\}=\{2,2,n\}\ (n\ge2),\quad \{2,3,3\},\quad \{2,3,4\}\quad \{2,3,5\},$$
in which case, the projective monodromy group $\bar G$ is isomorphic to the 
dihedral, tetrahedral, octahedral, and icosahedral group, of order
$$N=2n,\quad12,\quad24,\quad60,\qquad \frac2N=\frac1{k_0}+\frac1{k_1}+\frac1{k_\infty}-1,$$respectively.
In each case, the degrees of the basic invariants $\{P_0,P_1,P_\infty\}$ are
$$\{2,n,n\},\quad \{4,4,6\},\quad\{6,8,12\},\quad\{12,20,30\},$$
respectively. Since order of the $k$'s does matter, we name the triples $[k_0,k_1,k_\infty]$ as
$$\begin{array}{llll}
Dihedral: & D1=[2,2,n],&D2=[2,n,2],&D3=[n,2,2], \\[2mm]
Tetrahedral: & T1=[2,3,3],&T2=[3,2,3],&T3=[3,3,2], \\[2mm]
Octahedral:&O1=[2,3,4],&O2=[2,4,3],&O3=[3,2,4],\\[2mm]
&O4=[3,4,2],&O5=[4,2,3],&O6=[4,3,2], \\[2mm]
Icosahedral:&I1=[2,3,5],&I2=[2,3,5],&I3=[3,2,5],\\[2mm]
&I4=[3,5,2],&I5=[5,2,3],&I6=[5,3,2].\end{array}$$
For a triple $ [ k_0, k_1, k_\infty ]$, there are four triples $( a, b, c )$, up to the exchange of $a$ and $b$ (change of the sign of $\mu_\infty$):
$$\begin{array}{llll}
 \mathrm{i})~   &  \mu_0 > 0,~ \mu_1 > 0,~ \mu_\infty > 0,
&\quad\mathrm{ii})~  &  \mu_0 > 0,~ \mu_1 < 0,~ \mu_\infty < 0,\\[2mm]
 \mathrm{iii})~ &  \mu_0 < 0,~ \mu_1 > 0,~ \mu_\infty < 0,
&\quad\mathrm{iv})~  &  \mu_0 < 0,~ \mu_1 < 0,~ \mu_\infty > 0.\\
  \end{array}$$
\subsection{Theorem}
\begin{theorem}
Take a triple of real parameters $(a,b,c)$ satisfying
$$ |\mu_0|=\frac1{k_0},\quad |\mu_1|=\frac1{k_1},\quad |\mu_\infty|=\frac1{k_\infty},$$
where 
$$\mu_0=1-c,\quad \mu_1=c-a-b,\quad \mu_\infty=a-b,$$
and $[k_0,k_1,k_\infty]$ is a triple (listed above) of characteristics of a polyhedral group. The two solutions
$$ u(x) := F(a,b,c~;x),\quad v(x) := x^{1-c} F(a-c+1,b-c+1,2-c~;x),$$
of the hypergeometric differential equation $E(a,b,c)$ satisfy the identities
$$P_j(u(x),v(x))=f_j(x),\qquad j=0,1,\infty,$$
where each $P_j$ is a homogeneous polynomial (given in the next subsection) in $(u,v)$ of degree $N/k_j$ ($N$ is the order of the polyhedral group), and the functions $f_j=f_j(x)$ are given as follows:
$$\begin{array}{ccccc}
&\mbox{sign of $\mu_j$} & f_0&f_1&f_\infty\\[2mm]
\mathrm{i})& (+++)&
  \displaystyle{x^\frac{1}{k_0}}&(1-x)^\frac{1}{k_1}&1 \\[2mm]
\mathrm{ii})& (+--)&
\displaystyle{x^\frac{1}{k_0} (1-x)^{ - \frac{N}{k_0 k_1}}}&
(1-x)^{ \frac{1}{k_1} - \frac{N}{k_1 k_1}}
&(1-x)^{ - \frac{N}{k_1 k_\infty}}
\\[2mm]
\mathrm{iii})&(-+-)&
\displaystyle{ 
x^{ \frac{1}{k_0} - \frac{N}{k_0 k_0}}}&\displaystyle{ 
x^{ - \frac{N}{k_0 k_1}} (1-x)^\frac{1}{k_1}}&\displaystyle{ 
x^{ - \frac{N}{k_0 k_\infty}}}\\[2mm]
\mathrm{iv})& (--+)&
x^{ \frac{1}{k_0} - \frac{N}{k_0 k_0}} (1-x)^{ - \frac{N}{k_0 k_1}}&\quad
x^{ - \frac{N}{k_0 k_1}} (1-x)^{ \frac{1}{k_1} - \frac{N}{k_1 k_1}}&\quad
x^{ - \frac{N}{k_0 k_\infty}} (1-x)^{ - \frac{N}{k_1 k_\infty}}
\end{array}$$
\end{theorem}
\noindent
{\bf Sketch of Proof:} Once the theorem is formulated, we have only to substitute
$u$ and $v$ the power series solutions above in the polynomial $P_j(u,v)$, and to substitute the power function $(1-x)^*$ in the function $f_j$ its binomial expansion, and then to identify the coefficients. 
\begin{corollary} For each case, the image curve under the affine Schwarz map, and its inverse map are given as follows.
$$\begin{array}{cccc}
&\mbox{sign of $\mu_j$} &\quad\mbox{equation of the image curve}\quad&\mbox{inverse map}\\[2mm]
\mathrm{i})& (+++)&P_\infty - 1 = 0 &P_0^{k_0} = x \\[2mm]
\mathrm{ii})&(+--)& P_1^{\frac{N}{k_\infty}} - P_\infty^{\frac{N}{k_1} - 1} = 0 
 &P_0^{k_0} P_\infty^{-k_\infty} = x \\[2mm]
\mathrm{iii})&(-+-)& P_0^{\frac{N}{k_\infty}} - P_\infty^{\frac{N}{k_0} - 1} = 0
& P_0^{k_0} P_\infty^{-k_\infty} = x \\[2mm]
\mathrm{iv})&(--+)&\quad (P_0 P_1)^\frac{N}{k_\infty} - P_\infty^{\frac{N}{k_0} + \frac{N}{k_1} - 1} = 0
\quad&\quad P_0^{k_0} P_\infty^{-k_\infty} = x 
  \end{array}$$
\end{corollary}
\noindent
{\bf Proof:}
For each case, from the expressions in the theorem, it is obvious that the image curve $C$ satisfies the equation. We show that it is irreducible. The monodromy group $G$ of $asch$ is a cyclic extension of the monodromy group $\bar G$ of $sch$:
$$1\longrightarrow \Z/d\Z\longrightarrow G\longrightarrow \bar G\longrightarrow 1.$$ The index $d=[G,\bar G]$ is the sheet number of the covering $C\ni(u,v)\mapsto u:v\in Z$. This implies that the defining equation $F$ of $C$ can be written in the form $F=P-Q$, where $P$ and $Q$ are homogeneous in $(u,v)$ such that the degrees of the two polynomilas satisfy $\deg P-\deg Q=d$. Note that the equation posed in the Corollary is in this form with the degree-difference $N/k_\infty$. On the other hand, since one of the local exponents at $x=0$ and $x=1$ is 0, the sheet number $d$ is determined only by the exponents $\{a,b\}$ at $\infty$ as $d=N|a-b|$. Thus we can conclude that the equation posed is just the defining equation $F$, which is irreducible.  
\begin{remark} As polynomials in $(u,v)$, we have
$$P_0^{k_0}+P_1^{k_1}-P_\infty^{k_\infty}=0.$$
\end{remark}
\begin{remark}The result for case {\rm (i)} yields the remaining cases;
they are obtained by the use of the Kummer solutions around $x=1$ and the 
linear relation connecting these and the Kummer solutions $(u(x),v(x))$ around $x=0$.
\end{remark}
\noindent
{\bf Proof:}
Set
$$G(\lambda,\mu,\nu;x)=F(a,b,c;x),$$
where
$$\lambda=\mu_0=1-c,\quad\mu=\mu_1=c-a-b,\quad\nu=\mu_\infty=a-b.$$
Then the Kummer solutions around $x=0$ are
$$u=G(\lambda,\mu,\nu;x),\quad v=x^\lambda G(-\lambda,\mu,\nu;x),$$
and those at $x=1$ are
$$u_1=G(\mu,\lambda,\nu;1-x),\quad v_1=(1-x)^\mu G(-\mu,\lambda,\nu;1-x);$$
they are related as $(u,v)=(u_1,v_1)C,$ where $C=C(\lambda,\mu,\nu)$
$$=
\left(\begin{array}{cc}
\dfrac{\Gamma(1-\lambda)\Gamma(\mu)}
{\Gamma\left(\frac{1-\lambda+\mu-\nu}2\right)\Gamma\left(\frac{1-\lambda+\mu+\nu}2\right)}&
\dfrac{\Gamma(1+\lambda)\Gamma(\mu)}
{\Gamma\left(\frac{1+\lambda+\mu-\nu}2\right)\Gamma\left(\frac{1+\lambda+\mu+\nu}2\right)}\\[10mm]
\dfrac{\Gamma(1-\lambda)\Gamma(-\mu)}
{\Gamma\left(\frac{1-\lambda-\mu-\nu}2\right)\Gamma\left(\frac{1-\lambda-\mu+\nu}2\right)}&
\dfrac{\Gamma(1+\lambda)\Gamma(-\mu)}
{\Gamma\left(\frac{1+\lambda-\mu-\nu}2\right)\Gamma\left(\frac{1+\lambda-\mu+\nu}2\right)}
\end{array}\right);$$
note that the change of the sign of $\mu$ exchanges the rows.
\par\medskip\noindent
{\bf (i) $\Rightarrow$ (ii):} Set
$$\lambda'=\lambda,\quad \mu'=-\mu,\quad \nu'=\nu;$$
$u',v',u'_1,v'_1$ and $C'$ are defined using these parameters.
Then we have
$$u'_1=G(-\mu,\lambda,\nu;1-x)=(1-x)^{\mu'}v_1,\quad v'_1=(1-x)^{\mu'}G(\mu,\lambda,\nu;1-x)=(1-x)^{\mu'}u_1,$$
and so
$$(u'_1,v'_1)=(1-x)^{\mu'}(u_1,v_1)\left(\begin{array}{cc}0&1\\1&0\end{array}\right),$$and
$$\begin{array}{ll} (u',v')&=(u'_1,v'_1)C'=(1-x)^{\mu'}(u_1,v_1)\left(\begin{array}{cc}0&1\\1&0\end{array}\right)C'\\&=(1-x)^{\mu'}(u,v)C^{-1}\left(\begin{array}{cc}0&1\\1&0\end{array}\right)C'=(1-x)^{\mu'}(u,v).\end{array}$$
Thus we have
$$P_j(u',v')=(1-x)^{\mu'\deg P_j}P_j(u,v)=(1-x)^{\frac{-N}{k_1k_j}}f_j.$$

\par\medskip\noindent
{\bf (i) $ \Rightarrow $ (iii):} Set
$$\lambda'=-\lambda,\quad \mu'=\mu,\quad \nu'=\nu.$$
We have
$$u'=G(-\lambda,\mu,\nu)=x^{\lambda'}v,\quad v'=x^{\lambda'}u,$$
and so
$$P_j(v',u')=x^{\lambda'\deg P_j}P_j(u,v)=x^{\frac{-N}{k_0k_j}}f_j;$$
Note the order of $v'$ and $u'$ in $P_j(v',u')$.
\par\medskip\noindent
{\bf (i) $ \Rightarrow $ (iv):} It is a composition of the above two.

\subsection{Polynomials $ P_j $ and the explicit relations}
We tabulate, for each $[k_0,k_1,k_\infty]$, the polynomials $P_0,P_1,P_\infty$, and for the four possible parameter values $(a,b,c)$, the explicit form of the relations. Notational convention:
$$P_j=P_j(u,v),\quad P'_j=P_j(v,u).$$
Note that in the third and forth cases, $P'_j=P_j(v,u)$ are used.

\par\smallskip\noindent{$\bullet$\quad $D1=\left[ 2, 2, n \right] $}
$$\begin{array}{l}    P_0:= i \sum_{k=0}^{[\frac{n}{2}]} {n\choose2k+1} u^{n-2k-1} \left( \frac{1}{i} \frac{v}{n} \right)^{2k+1},\quad
  \ds P_1:= \sum_{k=0}^{[\frac{n}{2}]}{n\choose2k} u^{n-2k} \left( \frac{1}{i} \frac{v}{n} \right)^{2k}, \quad
  \ds P_\infty:= u^2 + \frac{1}{n^2} v^2.\end{array}$$
$$\begin{array}{llll}
 \left( \frac{1}{2n},-\frac{1}{2n}, \frac{1}{2} \right) & 
\left\{\begin{array}{ll}
  \ds P_0 &= x^\frac{1}{2} \\[2pt]
  \ds P_1  &= (1-x)^\frac{1}{2}\\[2pt]
  \ds P_\infty & = 1 \end{array}\right .
& \left( \frac{n-1}{2n}, \frac{n+1}{2n}, \frac{1}{2} \right) & 
\left\{\begin{array}{ll} 
  \ds P_0  &= x^\frac{1}{2} (1-x)^{-\frac{n}{2}} \\[2pt]
  \ds P_1  &= (1-x)^\frac{1-n}{2} \\[2pt]
  \ds P_\infty  &= (1-x)^{-1} 
       \end{array}\right .\\[8mm]
 \left( \frac{n-1}{2n}, \frac{n+1}{2n}, \frac{3}{2} \right) & 
 \left\{\begin{array}{ll}
  \ds P'_0&= x^\frac{1-n}{2} \\[2pt]
  \ds P'_1& = x^{-\frac{n}{2}} (1-x)^\frac{1}{2} \\[2pt]
  \ds P'_\infty&= x^{-1} \end{array}\right .
& \left(\frac{2n-1}{2n}, \frac{2n+1}{2n}, \frac{3}{2} \right) & 
   \left\{\begin{array}{ll}
 \ds P'_0&= x^\frac{1-n}{2} (1-x)^{-\frac{n}{2}} \\[2pt]
  \ds P'_1&= x^{-\frac{n}{2}} (1-x)^\frac{1-n}{2} \\[2pt]
  \ds P'_\infty&= x^{-1} (1-x)^{-1} 
 	  \end{array}\right .\end{array}$$
 \par\smallskip\noindent{$\bullet$\quad$D2= \left[ 2, n, 2 \right] $}
$$  \ds P_0:= \sum_{k=0}^{[\frac{n}{2}]} {n\choose 2k+1} u^{n-(2k+1)} \left( \frac{v}{n} \right)^{2k+1},\quad
  \ds P_1:= u^2 - \frac{1}{n^2} v^2,\quad
  \ds P_\infty:= \sum_{k=0}^{[\frac{n}{2}]} {n\choose 2k} u^{n-2k} \left( \frac{v}{n} \right)^{2k}.$$
$$
\begin{array}{llll}
 \left( \frac{n-1}{2n},-\frac{1}{2n}, \frac{1}{2} \right) &\left\{ 
 \begin{array}{ll}
  \ds P_0  &= x^\frac{1}{2} \\[2pt]
  \ds P_1  &= (1-x)^\frac{1}{n} \\[2pt]
  \ds P_\infty & = 1 
 \end{array} \right .
 &\left( \frac{1}{2n}, \frac{n+1}{2n}, \frac{1}{2} \right) &\left\{ 
 \begin{array}{ll}
  \ds P_0  &= x^\frac{1}{2} (1-x)^{-1} \\[2pt]
  \ds P_1  &= (1-x)^{-\frac{1}{n}} \\[2pt]
  \ds P_\infty  &= (1-x)^{-1} 
 \end{array}\right .\\[8mm]
 \left( \frac{n-1}{2n}, \frac{2n-1}{2n}, \frac{3}{2} \right) &\left\{ 
 \begin{array}{ll}
  \ds P'_0 &= x^\frac{1-n}{2} \\[2pt]
  \ds P'_1 &= x^{-1} (1-x)^\frac{1}{n} \\[2pt]
  \ds P'_\infty &= x^{-\frac{n}{2}} 
 \end{array}\right . &
 \left( \frac{2n+1}{2n}, \frac{n+1}{2n}, \frac{3}{2} \right) &\left\{ 
 \begin{array}{ll}
  \ds P'_0 &= x^\frac{1-n}{2} (1-x)^{-1} \\[2pt]
  \ds P'_1 &= x^{-1} (1-x)^{-\frac{1}{n}} \\[2pt]
  \ds P'_\infty &= x^{-\frac{n}{2}} (1-x)^{-1} 
 \end{array}\right .
\end{array}$$
\par\smallskip\noindent{$\bullet$\quad$D3= \left[ n, 2, 2 \right] $}
$$ \ds P_0  := u v,\quad  \ds P_1  := u^n - \frac{1}{4} v^n,\quad
  \ds P_\infty  := u^n + \frac{1}{4} v^n.$$
$$
\begin{array}{llll}
 \left( \frac{n-1}{2n}, -\frac{1}{2n}, \frac{n-1}{n} \right) &\left\{ 
 \begin{array}{ll}
  \ds P_0  &= x^\frac{1}{n} \\[2pt]
  \ds P_1  &= (1-x)^\frac{1}{2} \\[2pt]
  \ds P_\infty  &= 1 
 \end{array}  \right .
 &\left( \frac{n-1}{2n}, \frac{2n-1}{2n}, \frac{n-1}{n} \right) &\left\{ 
 \begin{array}{ll}
  \ds P_0  &= x^\frac{1}{n} (1-x)^{-1} \\[2pt]
  \ds P_1  &= (1-x)^\frac{1-n}{2} \\[2pt]
  \ds P_\infty  &= (1-x)^{-\frac{n}{2}} 
 \end{array} \right .\\[8mm]
 \left( \frac{1}{2n}, \frac{n+1}{2n}, \frac{n+1}{n} \right) &\left\{ 
 \begin{array}{ll}
  \ds P'_0 &= x^{-\frac{1}{n}} \\[2pt]
  \ds P'_1 &= x^{-1} (1-x)^\frac{1}{2} \\[2pt]
  \ds P'_\infty &= x^{-1} 
 \end{array}\right .  &
 \left( \frac{2n+1}{2n}, \frac{n+1}{2n}, \frac{n+1}{n} \right) &\left\{ 
 \begin{array}{ll}
  \ds P'_0 &= x^{-\frac{1}{n}} (1-x)^{-1} \\[2pt]
  \ds P'_1 &= x^{-1} (1-x)^\frac{1-n}{2} \\[2pt]
  \ds P'_\infty &= x^{-1} (1-x)^{-\frac{n}{2}} 
 \end{array} \right .
\end{array}$$
\par\smallskip\noindent{$\bullet$\quad$T1= \left[ 2, 3, 3 \right] $}
$$  \ds P_0  := u^5 v + \frac{1}{432} u v^5,\quad
  \ds P_1  := u^4 - \frac{1}{6} u^2 v^2 - \frac{1}{432} v^4,\quad
  \ds P_\infty  := u^4 + \frac{1}{6} u^2 v^2 - \frac{1}{432} v^4.$$
$$\begin{array}{llll}
 \left( \frac{1}{4},-\frac{1}{12}, \frac{1}{2} \right) &\left\{ 
 \begin{array}{llll}
  \ds P_0 &= x^\frac{1}{2} \\[2pt]
  \ds P_1 &= (1-x)^\frac{1}{3} \\[2pt]
  \ds P_\infty &= 1 
 \end{array}\right .  &
 \left( \frac{1}{4}, \frac{7}{12}, \frac{1}{2} \right) &\left\{ 
 \begin{array}{llll}
  \ds P_0 &= x^\frac{1}{2} (1-x)^{-2} \\[2pt]
  \ds P_1 &= (1-x)^{-1} \\[2pt]
  \ds P_\infty &= (1-x)^{ - \frac{4}{3}} 
 \end{array}\right . \\
 \left( \frac{5}{12}, \frac{3}{4}, \frac{3}{2} \right) &\left\{ 
 \begin{array}{llll}
  \ds P'_0&= x^{ - \frac{5}{2}} \\[2pt]
  \ds P'_1&= x^{-2} (1-x)^\frac{1}{3} \\[2pt]
  \ds P'_\infty&= x^{-2} 
 \end{array}\right .  &
 \left( \frac{13}{12}, \frac{3}{4}, \frac{3}{2} \right) &\left\{ 
 \begin{array}{llll}
  \ds P'_0&= x^{ - \frac{5}{2}} (1-x)^{-2} \\[2pt]
  \ds P'_1&= x^{-2} (1-x)^{-1} \\[2pt]
  \ds P'_\infty&= x^{-2} (1-x)^{ - \frac{4}{3}} 
 \end{array}\right . \\
\end{array}$$

\par\smallskip\noindent{$\bullet$\quad$T2=\left[ 3, 2, 3 \right] $}
$$
 \begin{array}{llll}
  \ds P_0  := u^3 v - \frac{1}{64} v^4,\quad
  \ds P_1  := u^6 - \frac{5}{16} u^3 v^3 - \frac{1}{512} v^6,\quad
  \ds P_\infty  := u^4 + \frac{1}{8} u v^3.
 \end{array}$$
$$\begin{array}{llll}
 \left( \frac{1}{4},-\frac{1}{12}, \frac{2}{3} \right) &\left\{ 
 \begin{array}{llll}
  \ds P_0 &= x^\frac{1}{3} \\[2pt]
  \ds P_1 &= (1-x)^\frac{1}{2} \\[2pt]
  \ds P_\infty &= 1 
 \end{array}\right .  &
 \left( \frac{5}{12}, \frac{3}{4}, \frac{2}{3} \right) &\left\{ 
 \begin{array}{llll}
  \ds P_0 &= x^\frac{1}{3} (1-x)^{-2} \\[2pt]
  \ds P_1 &= (1-x)^{ - \frac{5}{2}} \\[2pt]
  \ds P_\infty &= (1-x)^{-2} 
 \end{array}\right . \\
 \left( \frac{1}{4}, \frac{7}{12}, \frac{4}{3} \right) &\left\{ 
 \begin{array}{llll}
  \ds P'_0&= x^{-1} \\[2pt]
  \ds P'_1&= x^{-2} (1-x)^\frac{1}{2} \\[2pt]
  \ds P'_\infty&= x^{ - \frac{4}{3}} 
 \end{array}\right .  &
 \left( \frac{13}{12}, \frac{3}{4}, \frac{4}{3} \right) &\left\{ 
 \begin{array}{llll}
  \ds P'_0&= x^{-1} (1-x)^{-2} \\[2pt]
  \ds P'_1&= x^{-2} (1-x)^{ - \frac{5}{2}} \\[2pt]
  \ds P'_\infty&= x^{ - \frac{4}{3}} (1-x)^{-2} 
 \end{array}\right . \\
\end{array}$$
\par\smallskip\noindent{$\bullet$\quad$T3=\left[ 3, 3, 2 \right] $}
$$\begin{array}{llll}
  \ds P_0  := u^3 v + \frac{1}{64} v^4,\quad
  \ds P_1  := u^4 - \frac{1}{8} u v^3,\quad
  \ds P_\infty  := u^6 + \frac{5}{16} u^3 v^3 - \frac{1}{512} v^6 \\
 \end{array}$$
$$\begin{array}{llll}
 \left( \frac{5}{12},-\frac{1}{12}, \frac{2}{3} \right) &\left\{ 
 \begin{array}{llll}
  \ds P_0 &= x^\frac{1}{3} \\[2pt]
  \ds P_1 &= (1-x)^\frac{1}{3} \\[2pt]
  \ds P_\infty &= 1 
 \end{array}\right .  &
 \left( \frac{1}{4}, \frac{3}{4}, \frac{2}{3} \right) &\left\{ 
 \begin{array}{llll}
  \ds P_0 &= x^\frac{1}{3} (1-x)^{ - \frac{4}{3}} \\[2pt]
  \ds P_1 &= (1-x)^{-1} \\[2pt]
  \ds P_\infty &= (1-x)^{-2} 
 \end{array}\right . \\
 \left( \frac{1}{4}, \frac{3}{4}, \frac{4}{3} \right) &\left\{ 
 \begin{array}{llll}
  \ds P'_0&= x^{-1} \\[2pt]
  \ds P'_1&= x^{ - \frac{4}{3}} (1-x)^\frac{1}{3} \\[2pt]
  \ds P'_\infty&= x^{-2} 
 \end{array}\right .  &
 \left( \frac{13}{12}, \frac{7}{12}, \frac{4}{3} \right) &\left\{ 
 \begin{array}{llll}
  \ds P'_0&= x^{-1} (1-x)^{ - \frac{4}{3}} \\[2pt]
  \ds P'_1&= x^{ - \frac{4}{3}} (1-x)^{-1} \\[2pt]
  \ds P'_\infty&= x^{-2} (1-x)^{-2} 
 \end{array}\right . \\
\end{array}$$

\par\smallskip\noindent{$\bullet$\quad$O1= \left[ 2, 3, 4 \right] $}
$$ \begin{array}{l}
  \ds P_0  := u^{11} v - \frac{11}{432} u^9 v^3 + \frac{11}{3456} u^7 v^5 + \frac{11}{165888} u^5 v^7 - \frac{11}{47775744} u^3 v^9 + \frac{1}{254803968} u v^{11} \\[10pt]
  \ds P_1  := u^8 - \frac{7}{36} u^6 v^2 - \frac{7}{3456} u^4 v^4 - \frac{7}{82944} u^2 v^6 + \frac{1}{5308416} v^8 \\[10pt]
  \ds P_\infty  := u^6 + \frac{5}{48} u^4 v^2 - \frac{5}{2304} u^2 v^4 - \frac{1}{110592} v^6 \\
 \end{array}$$
$$\begin{array}{llll}
 \left( \frac{5}{24}, -\frac{1}{24}, \frac{1}{2} \right) &\left\{ 
 \begin{array}{llll}
  \ds P_0 &= x^\frac{1}{2} \\[2pt]
  \ds P_1 &= (1-x)^\frac{1}{3} \\[2pt]
  \ds P_\infty &= 1 
 \end{array}\right . &
 \left( \frac{7}{24}, \frac{13}{24}, \frac{1}{2} \right) &\left\{ 
 \begin{array}{llll}
  \ds P_0 &= x^\frac{1}{2} (1-x)^{-4} \\[2pt]
  \ds P_1 &= (1-x)^{ - \frac{7}{3}} \\[2pt]
  \ds P_\infty &= (1-x)^{-2} 
 \end{array}\right .\\
 \left( \frac{11}{24}, \frac{17}{24}, \frac{3}{2} \right) &\left\{ 
 \begin{array}{llll}
  \ds P'_0&= x^{ - \frac{11}{2}} \\[2pt]
  \ds P'_1&= x^{-4} (1-x)^\frac{1}{3} \\[2pt]
  \ds P'_\infty&= x^{-3} 
 \end{array}\right . &
 \left( \frac{25}{24}, \frac{19}{24}, \frac{3}{2} \right) &\left\{ 
 \begin{array}{llll}
  \ds P'_0&= x^{ - \frac{11}{2}} (1-x)^{-4} \\[2pt]
  \ds P'_1&= x^{-4} (1-x)^{ - \frac{7}{3}} \\[2pt]
  \ds P'_\infty&= x^{-3} (1-x)^{-2} 
 \end{array}\right .\\
\end{array}$$

\par\smallskip\noindent{$\bullet$\quad$O2=\left[ 2, 4, 3 \right] $}
$$ \begin{array}{l}
  \ds P_0  := u^{11} v + \frac{11}{432} u^9 v^3 + \frac{11}{3456} u^7 v^5 - \frac{11}{165888} u^5 v^7 - \frac{11}{47775744} u^3 v^9 - \frac{1}{254803968} u v^{11} \\[10pt]
  \ds P_1  := u^6 - \frac{5}{48} u^4 v^2 - \frac{5}{2304} u^2 v^4 + \frac{1}{110592} v^6 \\[10pt]
  \ds P_\infty  := u^8 + \frac{7}{36} u^6 v^2 - \frac{7}{3456} u^4 v^4 + \frac{7}{82944} u^2 v^6 + \frac{1}{5308416} v^8 \\
 \end{array}$$
$$\begin{array}{llll}
 \left( \frac{7}{24}, -\frac{1}{24}, \frac{1}{2} \right) &\left\{ 
 \begin{array}{llll}
  \ds P_0 &= x^\frac{1}{2} \\[2pt]
  \ds P_1 &= (1-x)^\frac{1}{4} \\[2pt]
  \ds P_\infty &= 1 
 \end{array}\right . &
 \left( \frac{5}{24}, \frac{13}{24}, \frac{1}{2} \right) &\left\{ 
 \begin{array}{llll}
  \ds P_0 &= x^\frac{1}{2} (1-x)^{-3} \\[2pt]
  \ds P_1 &= (1-x)^{ - \frac{5}{4}} \\[2pt]
  \ds P_\infty &= (1-x)^{-2} 
 \end{array}\right .\\
 \left( \frac{11}{24}, \frac{19}{24}, \frac{3}{2} \right) &\left\{ 
 \begin{array}{llll}
  \ds P'_0&= x^{ - \frac{11}{2}} \\[2pt]
  \ds P'_1&= x^{-3} (1-x)^\frac{1}{4} \\[2pt]
  \ds P'_\infty&= x^{-4} 
 \end{array}\right . &
 \left( \frac{25}{24}, \frac{17}{24}, \frac{3}{2} \right) &\left\{ 
 \begin{array}{llll}
  \ds P'_0&= x^{ - \frac{11}{2}} (1-x)^{-3} \\[2pt]
  \ds P'_1&= x^{-3} (1-x)^{ - \frac{5}{4}} \\[2pt]
  \ds P'_\infty&= x^{-4} (1-x)^{-4} 
 \end{array}\right .\\
\end{array}$$

\par\smallskip\noindent{$\bullet$\quad$ O3=\left[ 3, 2, 4 \right] $}
$$\begin{array}{llll}
  \ds P_0  := u^7 v - \frac{7}{256} u^4 v^4 - \frac{1}{8192} u v^7,\quad
  \ds P_1  := u^{12} - \frac{11}{32} u^9 v^3 - \frac{11}{262144} u^3 v^9 - \frac{1}{67108864} v^{12} \\[10pt]
  \ds P_\infty  := u^6 + \frac{5}{64} u^3 v^3 - \frac{1}{8192} v^6 \\
 \end{array}$$
$$\begin{array}{llll}
 \left( \frac{5}{24}, -\frac{1}{24}, \frac{2}{3} \right) &\left\{ 
 \begin{array}{llll}
  \ds P_0 &= x^\frac{1}{3} \\[2pt]
  \ds P_1 &= (1-x)^\frac{1}{2} \\[2pt]
  \ds P_\infty &= 1 
 \end{array}\right . &
 \left( \frac{11}{24}, \frac{17}{24}, \frac{2}{3} \right) &\left\{ 
 \begin{array}{llll}
  \ds P_0 &= x^\frac{1}{3} (1-x)^{-4} \\[2pt]
  \ds P_1 &= (1-x)^{ - \frac{11}{2}} \\[2pt]
  \ds P_\infty &= (1-x)^{-3} 
 \end{array}\right .\\
 \left( \frac{7}{24}, \frac{13}{24}, \frac{4}{3} \right) &\left\{ 
 \begin{array}{llll}
  \ds P'_0&= x^{ - \frac{7}{3}} \\[2pt]
  \ds P'_1&= x^{-4} (1-x)^\frac{1}{2} \\[2pt]
  \ds P'_\infty&= x^{-2} 
 \end{array}\right . &
 \left( \frac{25}{24}, \frac{19}{24}, \frac{4}{3} \right) &\left\{ 
 \begin{array}{llll}
  \ds P'_0&= x^{ - \frac{7}{3}} (1-x)^{-4} \\[2pt]
  \ds P'_1&= x^{-4} (1-x)^{ - \frac{11}{2}} \\[2pt]
  \ds P'_\infty&= x^{-2} (1-x)^{-3} 
 \end{array}\right .\\
\end{array}$$

\par\smallskip\noindent{$\bullet$\quad$O4= \left[ 3, 4, 2 \right] $}
$$ \begin{array}{llll}
  \ds P_0  := u^7 v + \frac{7}{256} u^4 v^4 - \frac{1}{8192} u v^7,\quad
  \ds P_1  := u^6 - \frac{5}{64} u^3 v^3 - \frac{1}{8192} v^6 \\[10pt]
  \ds P_\infty  := u^{12} + \frac{11}{32} u^9 v^3 + \frac{11}{262144} u^3 v^9 - \frac{1}{67108864} v^{12} \\
 \end{array}$$
$$\begin{array}{llll}
 \left( \frac{11}{24}, -\frac{1}{24}, \frac{2}{3} \right) &\left\{ 
 \begin{array}{llll}
  \ds P_0 &= x^\frac{1}{3} \\[2pt]
  \ds P_1 &= (1-x)^\frac{1}{4} \\[2pt]
  \ds P_\infty &= 1 
 \end{array}\right . &
 \left( \frac{5}{24}, \frac{17}{24}, \frac{2}{3} \right) &\left\{ 
 \begin{array}{llll}
  \ds P_0 &= x^\frac{1}{3} (1-x)^{-2} \\[2pt]
  \ds P_1 &= (1-x)^{ - \frac{5}{4}} \\[2pt]
  \ds P_\infty &= (1-x)^{-3} 
 \end{array}\right .\\
 \left( \frac{7}{24}, \frac{19}{24}, \frac{4}{3} \right) &\left\{ 
 \begin{array}{llll}
  \ds P'_0&= x^{ - \frac{7}{3}} \\[2pt]
  \ds P'_1&= x^{-2} (1-x)^\frac{1}{4} \\[2pt]
  \ds P'_\infty&= x^{-4} 
 \end{array}\right . &
 \left( \frac{25}{24}, \frac{13}{24}, \frac{4}{3} \right) &\left\{ 
 \begin{array}{llll}
  \ds P'_0&= x^{ - \frac{7}{3}} (1-x)^{-2} \\[2pt]
  \ds P'_1&= x^{-2} (1-x)^{ - \frac{5}{4}} \\[2pt]
  \ds P'_\infty&= x^{-4} (1-x)^{-3} 
 \end{array}\right .\\
\end{array}$$

\par\smallskip\noindent{$\bullet$\quad$O5= \left[ 4, 2, 3 \right] $}
$$ \begin{array}{llll}
  \ds P_0  := u^5 v - \frac{1}{108} u v^5,\quad
  \ds P_1  := u^{12} - \frac{11}{36} u^8 v^4 - \frac{11}{3888} u^4 v^8 + \frac{1}{1259712} v^{12} \\[10pt]
  \ds P_\infty  := u^8 + \frac{7}{54} u^4 v^4 + \frac{1}{11664} v^8 \\
 \end{array}$$
$$\begin{array}{llll}
 \left( \frac{7}{24}, -\frac{1}{24}, \frac{3}{4} \right) &\left\{ 
 \begin{array}{llll}
  \ds P_0 &= x^\frac{1}{4} \\[2pt]
  \ds P_1 &= (1-x)^\frac{1}{2} \\[2pt]
  \ds P_\infty &= 1 
 \end{array}\right . &
 \left( \frac{11}{24}, \frac{19}{24}, \frac{3}{4} \right) &\left\{ 
 \begin{array}{llll}
  \ds P_0 &= x^\frac{1}{2} (1-x)^{-3} \\[2pt]
  \ds P_1 &= (1-x)^{ - \frac{11}{2}} \\[2pt]
  \ds P_\infty &= (1-x)^{-4} 
 \end{array}\right .\\
 \left( \frac{5}{24}, \frac{13}{24}, \frac{5}{4} \right) &\left\{ 
 \begin{array}{llll}
  \ds P'_0&= x^{ - \frac{5}{4}} \\[2pt]
  \ds P'_1&= x^{-3} (1-x)^\frac{1}{2} \\[2pt]
  \ds P'_\infty&= x^{-2} 
 \end{array}\right . &
 \left( \frac{25}{24}, \frac{17}{24}, \frac{5}{4} \right) &\left\{ 
 \begin{array}{llll}
  \ds P'_0&= x^{ - \frac{5}{4}} (1-x)^{-3} \\[2pt]
  \ds P'_1&= x^{-3} (1-x)^{ - \frac{11}{2}} \\[2pt]
  \ds P'_\infty&= x^{-2} (1-x)^{-4} 
 \end{array}\right .\\
\end{array}$$

\par\smallskip\noindent{$\bullet$\quad$O6= \left[ 4, 3, 2 \right] $}
$$\begin{array}{llll}
  \ds P_0  := u^5 v + \frac{1}{108} u v^5,\quad
  \ds P_1  := u^8 - \frac{7}{54} u^4 v^4 + \frac{1}{11664} v^8 \\[10pt]
  \ds P_\infty  := u^{12} + \frac{11}{36} u^8 v^4 - \frac{11}{3888} u^4 v^8 - \frac{1}{1259712} v^{12} \\
 \end{array}$$
$$\begin{array}{llll}
 \left( \frac{11}{24}, -\frac{1}{24}, \frac{3}{4} \right) &\left\{ 
 \begin{array}{llll}
  \ds P_0 &= x^\frac{1}{4} \\[2pt]
  \ds P_1 &= (1-x)^\frac{1}{3} \\[2pt]
  \ds P_\infty &= 1 
 \end{array}\right . &
 \left( \frac{7}{24}, \frac{19}{24}, \frac{3}{4} \right) &\left\{ 
 \begin{array}{llll}
  \ds P_0 &= x^\frac{1}{4} (1-x)^{-2} \\[2pt]
  \ds P_1 &= (1-x)^{ - \frac{7}{3}} \\[2pt]
  \ds P_\infty &= (1-x)^{-4} 
 \end{array}\right .\\
 \left( \frac{5}{24}, \frac{17}{24}, \frac{5}{4} \right) &\left\{ 
 \begin{array}{llll}
  \ds P'_0&= x^{ - \frac{5}{4}} \\[2pt]
  \ds P'_1&= x^{-2} (1-x)^\frac{1}{3} \\[2pt]
  \ds P'_\infty&= x^{-3} 
 \end{array}\right . &
 \left( \frac{25}{24}, \frac{13}{24}, \frac{5}{4} \right) &\left\{ 
 \begin{array}{llll}
  \ds P'_0&= x^{ - \frac{5}{4}} (1-x)^{-2} \\[2pt]
  \ds P'_1&= x^{-2} (1-x)^{ - \frac{7}{3}} \\[2pt]
  \ds P'_\infty&= x^{-3} (1-x)^{-4} 
 \end{array}\right .\\
\end{array}$$

\par\smallskip\noindent{$\bullet$\quad$I1= \left[ 2, 3, 5 \right] $}
$$ \begin{array}{l}
  \ds P_0  := u^{29} v - \frac{29}{675} u^{27} v^3 + \frac{1769}{450000} u^{25} v^5 + \frac{29}{337500} u^{23} v^7 + \frac{667}{540000000} u^{21} v^9 - \frac{667}{72900000000} u^{19} v^{11} \\[10pt]
  \ds~~~~~~~~~~ + \frac{12673}{29160000000000} u^{17} v^{13} - \frac{12673}{524880000000000000} v^{13} u^{17} + \frac{667}{23619600000000000000} u^{11} v^{19} \\[10pt]
  \ds~~~~~~~~~~ - \frac{667}{3149280000000000000000} u^9 v^{21} - \frac{29}{35429400000000000000000} u^7 v^{23} \\[10pt]
  \ds~~~~~~~~~~ - \frac{1769}{850305600000000000000000000} u^5 v^{25} + \frac{29}{22958251200000000000000000000} u^3 v^{27} \\[10pt]
  \ds~~~~~~~~~~ - \frac{1}{612220032000000000000000000000} u v^{29} \\[10pt]
  \ds P_1  := u^{20} - \frac{19}{90} u^{18} v^2 - \frac{19}{18000} u^{16} v^4 - \frac{19}{135000} u^{14} v^6 - \frac{247}{162000000} u^{12} v^8 + \frac{247}{14580000000} u^{10} v^{10} \\[10pt]
  \ds~~~~~~~~~~ - \frac{247}{2916000000000} u^8 v^{12} - \frac{19}{43740000000000} u^6 v^{14} - \frac{19}{104976000000000000} u^4 v^{16} \\[10pt]
  \ds~~~~~~~~~~ - \frac{19}{9447840000000000000} u^2 v^{18} + \frac{1}{1889568000000000000000} v^{20},\\[10pt]
  \ds P_\infty  := u^{12} + \frac{11}{150} u^{10} v^2 - \frac{11}{6000} u^8 v^4 - \frac{11}{1350000} u^6 v^6 - \frac{11}{108000000} u^4 v^8 \\[10pt]
  \ds~~~~~~~~~~~~~ + \frac{11}{48600000000} u^2 v^{10} + \frac{1}{5832000000000} v^{12}.\\
 \end{array}$$
$$\begin{array}{llll}
 \left( \frac{11}{60}, -\frac{1}{60}, \frac{1}{2} \right) &\left\{ 
 \begin{array}{llll}
  \ds P_0 &= x^\frac{1}{2} \\[2pt]
  \ds P_1 &= (1-x)^\frac{1}{3} \\[2pt]
  \ds P_\infty &= 1 
 \end{array}\right . &
 \left( \frac{19}{60}, \frac{31}{60}, \frac{1}{2} \right) &\left\{ 
 \begin{array}{llll}
  \ds P_0 &= x^\frac{1}{2} (1-x)^{-10} \\[2pt]
  \ds P_1 &= (1-x)^{ - \frac{19}{3}} \\[2pt]
  \ds P_\infty &= (1-x)^{-4} 
 \end{array}\right .\\
 \left( \frac{29}{60}, \frac{41}{60}, \frac{3}{2} \right) &\left\{ 
 \begin{array}{llll}
  \ds P'_0&= x^{ - \frac{29}{2}} \\[2pt]
  \ds P'_1&= x^{-10} (1-x)^\frac{1}{3} \\[2pt]
  \ds P'_\infty&= x^{-6} 
 \end{array}\right . &
 \left( \frac{61}{60}, \frac{49}{60}, \frac{3}{2} \right) &\left\{ 
 \begin{array}{llll}
  \ds P'_0&= x^{ - \frac{29}{2}} (1-x)^{-10} \\[2pt]
  \ds P'_1&= x^{-10} (1-x)^{ - \frac{19}{3}} \\[2pt]
  \ds P'_\infty&= x^{-6} (1-x)^{-4} 
 \end{array}\right .\\
\end{array}$$

\par\smallskip\noindent{$\bullet$\quad$I6=\left[ 5, 3, 2 \right] $}
$$ \begin{array}{llll}
  \ds P_0  := u^{11} v + \frac{11}{1728} u^6 v^6 - \frac{1}{2985984} u v^{11} \\[10pt]
  \ds P_1  := u^{20} - \frac{19}{144} u^{15} v^5 + \frac{247}{1492992} u^{10} v^{10} + \frac{19}{429981696} u^5 v^{15} + \frac{1}{8916100448256} v^{20} \\[10pt]
  \ds P_\infty  := u^{30} + \frac{29}{96} u^{25} v^5 - \frac{3335}{995328} u^{20} v^{10} - \frac{3335}{2972033482752} u^{10} v^{20} - \frac{29}{855945643032576} u^5 v^{25} \\[10pt]
  \ds~~~~~~~~~~~~~ + \frac{1}{26623333280885243904} v^{30} \\
 \end{array}$$
$$\begin{array}{llll}
 \left( \frac{29}{60}, -\frac{1}{60}, \frac{4}{5} \right) &\left\{ 
 \begin{array}{llll}
  \ds P_0 &= x^\frac{1}{5} \\[2pt]
  \ds P_1 &= (1-x)^\frac{1}{3} \\[2pt]
  \ds P_\infty &= 1 
 \end{array}\right . &
 \left( \frac{19}{60}, \frac{49}{60}, \frac{4}{5} \right) &\left\{ 
 \begin{array}{llll}
  \ds P_0 &= x^\frac{1}{5} (1-x)^{-4} \\[2pt]
  \ds P_1 &= (1-x)^{ - \frac{19}{3}} \\[2pt]
  \ds P_\infty &= (1-x)^{-10} 
 \end{array}\right .\\
 \left( \frac{11}{60}, \frac{41}{60}, \frac{6}{5} \right) &\left\{ 
 \begin{array}{llll}
  \ds P'_0&= x^{ - \frac{11}{5}} \\[2pt]
  \ds P'_1&= x^{-4} (1-x)^\frac{1}{3} \\[2pt]
  \ds P'_\infty&= x^{-6} 
 \end{array}\right . &
 \left( \frac{61}{60}, \frac{31}{60}, \frac{6}{5} \right) &\left\{ 
 \begin{array}{llll}
  \ds P'_0&= x^{ - \frac{11}{5}} (1-x)^{-4} \\[2pt]
  \ds P'_1&= x^{-4} (1-x)^{ - \frac{19}{3}} \\[2pt]
  \ds P'_\infty&= x^{-6} (1-x)^{-10} 
 \end{array}\right .\\
\end{array}$$
For the following cases, consult \cite{Kob}:
\par\smallskip\noindent{$I2= \left[ 2, 5, 3 \right] $}, {$I3= \left[ 3, 2, 5 \right] $}, {$I4= \left[ 3, 5, 2 \right] $}, {$I5= \left[ 5, 2, 3 \right] $}.

\subsection{Some other cases}
From the Schwarz' list of the parameters (mod $1$) of the hypergeometric equation with finite monodromy group, we pick up some of them. Note that in the identity of the form $$P(u(x),v(x))=(\mbox{a rational function in $x$})\cdot x^*\cdot(1-x)^*,$$ the rational function is not necessarily 1 any more.
\par\medskip\noindent{$\bullet$\quad$T4=\left[ 3, 3, \frac{3}{2} \right]$}
$$  \ds P_0  := u^3 v - \frac{1}{16} v^4,\quad
  \ds P_1  := u^6 - \frac{5}{4} u^3 v^3 - \frac{1}{32} v^6,\quad
  \ds P_\infty  := u^4 + \frac{1}{2} u v^3.$$
$$\begin{array}{llll}
 \left( \frac{1}{2},-\frac{1}{6}, \frac{2}{3} \right) &\left\{ 
 \begin{array}{llll}
  \ds P_0 &= x^\frac{1}{3} (1-x)^\frac{1}{3} \\[2pt]
  \ds P_1 &= 1 - 2 x \\[2pt]
  \ds P_\infty &= 1 
 \end{array}\right . &
 \left( \frac{1}{6}, \frac{5}{6}, \frac{2}{3} \right) &\left\{ 
 \begin{array}{llll}
  \ds P_0 &= x^\frac{1}{3} (1-x)^{-1} \\[2pt]
  \ds P_1 &= (1-x)^{-2} (1 - 2 x) \\[2pt]
  \ds P_\infty &= (1-x)^{ - \frac{4}{3}} 
 \end{array}\right .\\
 \left( \frac{1}{6}, \frac{5}{6}, \frac{4}{3} \right) &\left\{ 
 \begin{array}{llll}
  \ds P'_0&= x^{-1} (1-x)^\frac{1}{3} \\[2pt]
  \ds P'_1&= x^{-2} (1 - 2 x) \\[2pt]
  \ds P'_\infty&= x^{ - \frac{4}{3}} 
 \end{array}\right . &
 \left( \frac{7}{6}, \frac{1}{2}, \frac{4}{3} \right) &\left\{ 
 \begin{array}{llll}
  \ds P'_0&= x^{-1} (1-x)^{-1} \\[2pt]
  \ds P'_1&= x^{-2} (1-x)^{-2} (1 - 2 x) \\[2pt]
  \ds P'_\infty&= x^{ - \frac{4}{3}} (1-x)^{ - \frac{4}{3}} 
 \end{array}\right .\\
\end{array}$$

\par\smallskip\noindent{$\bullet$\quad$T5=\left[ 3, \frac{3}{2}, 3 \right] $}
$$\begin{array}{l}
  \ds P_0  := u^3 v + \frac{1}{16} v^4,\quad
  \ds P_1  := u^4 - \frac{1}{2} u v^3,\quad
  \ds P_\infty  := u^6 + \frac{5}{4} u^3 v^3 - \frac{1}{32} v^6.
 \end{array}$$
$$\begin{array}{llll}
 \left( \frac{1}{6},-\frac{1}{6}, \frac{2}{3} \right) &\left\{ 
 \begin{array}{llll}
  \ds P_0 &= x^\frac{1}{3} \\[2pt]
  \ds P_1 &= (1-x)^\frac{2}{3} \\[2pt]
  \ds P_\infty &= 1 + x 
 \end{array}\right . &
 \left( \frac{1}{2}, \frac{5}{6}, \frac{2}{3} \right) &\left\{ 
 \begin{array}{llll}
  \ds P_0 &= x^\frac{1}{3} (1-x)^{-\frac{8}{3}} \\[2pt]
  \ds P_1 &= (1-x)^{-2} \\[2pt]
  \ds P_\infty &= (1-x)^{-4} (1 + x) 
 \end{array}\right .\\
 \left( \frac{7}{6}, \frac{5}{6}, \frac{4}{3} \right) &\left\{ 
 \begin{array}{llll}
  \ds P'_0&= x^{-1} \\[2pt]
  \ds P'_1&= x^{ - \frac{4}{3}} (1-x)^\frac{2}{3} \\[2pt]
  \ds P'_\infty&= x^{-2} ( 1 + x ) 
 \end{array}\right . &
 \left( \frac{7}{6}, \frac{5}{6}, \frac{4}{3} \right) &\left\{ 
 \begin{array}{llll}
  \ds P'_0&= x^{-1} (1-x)^{ - \frac{8}{3}} \\[2pt]
  \ds P'_1&= x^{ - \frac{4}{3}} (1-x)^{-2} \\[2pt]
  \ds P'_\infty&= x^{-2} (1-x)^{-4} (1 + x) 
 \end{array}\right .\\
\end{array}$$

\par\smallskip\noindent{$\bullet$\quad$T6= \left[ \frac{3}{2}, 3, 3 \right] $}
$$ \begin{array}{l}
  \ds P_0  := u^3 v + \frac{1}{256} v^4,\quad
  \ds P_1  := u^4 - \frac{1}{32} u v^3,\quad
  \ds P_\infty  := u^6 + \frac{5}{64} u^3 v^3 - \frac{1}{8192} v^6.
 \end{array}$$
$$\begin{array}{llll}
 \left( \frac{1}{6},-\frac{1}{6}, \frac{1}{3} \right) &\left\{ 
 \begin{array}{llll}
  \ds P_0 &= x^\frac{2}{3} \\[2pt]
  \ds P_1 &= (1-x)^\frac{1}{3} \\[2pt]
  \ds P_\infty &= 1 - \frac{1}{2} x 
 \end{array}\right . &
 \left( \frac{1}{6}, \frac{1}{2}, \frac{1}{3} \right) &\left\{ 
 \begin{array}{llll}
  \ds P_0 &= x^\frac{2}{3} (1-x)^{-\frac{4}{3}} \\[2pt]
  \ds P_1 &= (1-x)^{-1} \\[2pt]
  \ds P_\infty &= (1-x)^{-2} ( 1 - \frac{x}{2} ) 
 \end{array}\right .\\
 \left( \frac{1}{2}, \frac{5}{6}, \frac{5}{3} \right) &\left\{ 
 \begin{array}{llll}
  \ds P'_0&= x^{-2} \\[2pt]
  \ds P'_1&= x^{ - \frac{8}{3}} (1-x)^\frac{1}{3} \\[2pt]
  \ds P'_\infty&= x^{-4} ( 1 - \frac{x}{2} ) 
 \end{array}\right . &
 \left( \frac{7}{6}, \frac{5}{6}, \frac{5}{3} \right) &\left\{ 
 \begin{array}{llll}
  \ds P'_0&= x^{-2} (1-x)^{ - \frac{4}{3}} \\[2pt]
  \ds P'_1&= x^{ - \frac{8}{3}} (1-x)^{-1} \\[2pt]
  \ds P'_\infty&= x^{-4} (1-x)^{-2} ( 1 - \frac{x}{2} ) 
 \end{array}\right .\\
\end{array}$$

\par\smallskip\noindent{$\bullet$\quad$O7= \ds \left[ 4, 4, \frac{3}{2} \right] $}
$$ \begin{array}{l}
  \ds P_0  := u^5 v - \frac{1}{27} u v^5,\quad
  \ds P_1  := u^{12} - \frac{11}{9} u^8 v^4 - \frac{11}{243} u^4 v^8 + \frac{1}{19683} v^{12},\quad
  \ds P_\infty  := u^8 + \frac{14}{27} u^4 v^4 + \frac{1}{729} v^8 \\
 \end{array}$$
$$\begin{array}{llll}
 \left( \frac{7}{12}, -\frac{1}{12}, \frac{3}{4} \right) &\left\{ 
 \begin{array}{llll}
  \ds P_0 &= x^\frac{1}{4}(1-x)^\frac{1}{4} \\[2pt]
  \ds P_1 &= 1 - 2x \\[2pt]
  \ds P_\infty &= 1 
 \end{array}\right . &
 \left( \frac{1}{6}, \frac{5}{6}, \frac{3}{4} \right) &\left\{ 
 \begin{array}{llll}
  \ds P_0 &= x^\frac{1}{4} (1-x)^{ - \frac{5}{4}} \\[2pt]
  \ds P_1 &= (1-x)^{-3} ( 1 - 2x ) \\[2pt]
  \ds P_\infty &= (1-x)^{-2} 
 \end{array}\right .\\
 \left( \frac{1}{6}, \frac{5}{6}, \frac{5}{4} \right) &\left\{ 
 \begin{array}{llll}
  \ds P'_0&= x^{ - \frac{5}{4}} (1-x)^\frac{1}{4} \\[2pt]
  \ds P'_1&= x^{-3} ( 1 - 2x ) \\[2pt]
  \ds P'_\infty&= x^{-2} 
 \end{array}\right . &
 \left( \frac{13}{12}, \frac{5}{12}, \frac{5}{4} \right) &\left\{ 
 \begin{array}{llll}
  \ds P'_0&= x^{ - \frac{5}{4}} (1-x)^{ - \frac{5}{4}} \\[2pt]
  \ds P'_1&= x^{-3} (1-x)^{-3} ( 1 - 2x ) \\[2pt]
  \ds P'_\infty&= x^{-2} (1-x)^{-2} ( 1 - \frac{x}{2} ) 
 \end{array}\right .\\
\end{array}$$

\par\smallskip\noindent{$\bullet$\quad$O8= \left[ 4, \frac{3}{2}, 4 \right] $}
$$ \begin{array}{l}
  \ds P_0  := u^5 v + \frac{1}{27} u v^5,\quad
  \ds P_1  := u^{12} + \frac{11}{9} u^8 v^4 - \frac{11}{243} u^4 v^8 - \frac{1}{19683} v^{12},\quad
  \ds P_\infty  := u^8 - \frac{14}{27} u^4 v^4 + \frac{1}{729} v^8.
 \end{array}$$
$$\begin{array}{llll}
 \left( \frac{1}{6}, -\frac{1}{12}, \frac{3}{4} \right) &\left\{ 
 \begin{array}{llll}
  \ds P_0 &= x^\frac{1}{4} \\[2pt]
  \ds P_1 &= (1-x)^\frac{2}{3} \\[2pt]
  \ds P_\infty &= 1 + x 
 \end{array}\right . &
 \left( \frac{7}{12}, \frac{5}{6}, \frac{3}{4} \right) &\left\{ 
 \begin{array}{llll}
  \ds P_0 &= x^\frac{1}{4} (1-x)^{-4} \\[2pt]
  \ds P_1 &= (1-x)^{ - \frac{14}{3}} \\[2pt]
  \ds P_\infty &= (1-x)^{-8} (1 + x) 
 \end{array}\right .\\
 \left( \frac{1}{6}, \frac{5}{12}, \frac{5}{4} \right) &\left\{ 
 \begin{array}{llll}
  \ds P'_0&= x^{ - \frac{5}{4}} \\[2pt]
  \ds P'_1&= x^{-2} (1-x)^\frac{2}{3} \\[2pt]
  \ds P'_\infty&= x^{-3} ( 1 + x ) 
 \end{array}\right . &
 \left( \frac{13}{12}, \frac{5}{6}, \frac{5}{4} \right) &\left\{ 
 \begin{array}{llll}
  \ds P'_0&= x^{ - \frac{5}{4}} (1-x)^{-4} \\[2pt]
  \ds P'_1&= x^{-2} (1-x)^{ - \frac{14}{3}} \\[2pt]
  \ds P'_\infty&= x^{-3} (1-x)^{-8} ( 1 + x ) 
 \end{array}\right .\\
\end{array}$$

\par\smallskip\noindent{$\bullet$\quad$O9=\left[ \frac{3}{2}, 4, 4 \right] $}
$$ \begin{array}{l}
  \ds P_0  := u^7 v + \frac{7}{1024} u^4 v^4 - \frac{1}{131072} u v^7,\quad
  \ds P_1  := u^6 - \frac{5}{256} u^3 v^3 - \frac{1}{131072} v^6 \\[10pt]
  \ds P_\infty  := u^{12} + \frac{11}{128} u^9 v^3 + \frac{11}{16777216} u^3 v^9 - \frac{1}{17179869184} v^{12} \\
 \end{array}$$
$$\begin{array}{llll}
 \left( \frac{1}{6}, -\frac{1}{12}, \frac{1}{3} \right) &\left\{ 
 \begin{array}{llll}
  \ds P_0 &= x^\frac{2}{3} \\[2pt]
  \ds P_1 &= (1-x)^\frac{1}{4} \\[2pt]
  \ds P_\infty &= 1 - \frac{1}{2} x 
 \end{array}\right . &
 \left( \frac{1}{6}, \frac{5}{12}, \frac{1}{3} \right) &\left\{ 
 \begin{array}{llll}
  \ds P_0 &= x^\frac{2}{3} (1-x)^{-2} \\[2pt]
  \ds P_1 &= (1-x)^{ - \frac{5}{4}} \\[2pt]
  \ds P_\infty &= (1-x)^{-3} (1 - \frac{1}{2} x) 
 \end{array}\right .\\
 \left( \frac{7}{12}, \frac{5}{6}, \frac{5}{3} \right) &\left\{ 
 \begin{array}{llll}
  \ds P'_0&= x^{ - \frac{14}{3}} \\[2pt]
  \ds P'_1&= x^{-4} (1-x)^\frac{1}{4} \\[2pt]
  \ds P'_\infty&= x^{-8} (1 - \frac{1}{2} x) 
 \end{array}\right . &
 \left( \frac{13}{12}, \frac{5}{6}, \frac{5}{3} \right) &\left\{ 
 \begin{array}{llll}
  \ds P'_0&= x^{ - \frac{14}{3}} (1-x)^{-2} \\[2pt]
  \ds P'_1&= x^{-4} (1-x)^{ - \frac{5}{4}} \\[2pt]
  \ds P'_\infty&= x^{-8} (1-x)^{-3} (1 - \frac{1}{2} x) 
 \end{array}\right .
\end{array}$$

 \par\smallskip\noindent{$\bullet$\quad$I7= \left[ 3, 3, \frac{5}{2} \right] $}
$$ \begin{array}{llll}
  \ds P_0 ( u, v ) := u^{19} v - \frac{57}{400} u^{16} v^4 - \frac{57}{25000} u^{13} v^7 - \frac{247}{2500000} u^{10} v^{10} + \frac{57}{312500000} u^7 v^{13} \\[10pt]
  \ds~~~~~~~~~~~ - \frac{57}{62500000000} u^4 v^{16} - \frac{1}{1953125000000} u v^{19} \\[10pt]
  \ds P_1 ( u, v ) := u^{30} - \frac{29}{20} u^{27} v^3 + \frac{783}{20000} u^{24} v^6 - \frac{2001}{500000} u^{21} v^9 - \frac{38019}{250000000} u^{18} v^{12} \\[10pt]
  \ds~~~~~~~~~~~ - \frac{38019}{3125000000000} u^{12} v^{18} + \frac{2001}{78125000000000} u^9 v^{21} + \frac{783}{39062500000000000} u^6 v^{24} \\[10pt]
  \ds~~~~~~~~~~~ + \frac{29}{488281250000000000} u^3 v^{27} + \frac{1}{305175781250000000000} v^{30} \\[10pt]
  \ds P_\infty ( u, v ) := u^{12} + \frac{11}{50} u^9 v^3 - \frac{3}{12500} u^6 v^6 - \frac{11}{625000} u^3 v^9 + \frac{1}{156250000} v^{12} \\
 \end{array}$$
$$\begin{array}{llll}
 \left( \frac{11}{30}, -\frac{1}{30}, \frac{2}{3} \right) &\left\{ 
 \begin{array}{llll}
  \ds P_0 ( u, v ) = x^\frac{1}{3} (1-x)^\frac{1}{3} \\[2pt]
  \ds P_1 ( u, v ) = 1-2x \\[2pt]
  \ds P_\infty ( u, v ) = 1 
 \end{array}\right . &
 \left( \frac{7}{10}, \frac{3}{10}, \frac{2}{3} \right) &\left\{ 
 \begin{array}{llll}
  \ds P_0 ( u, v ) = x^\frac{1}{3} (1-x)^{-\frac{19}{3}} \\[2pt]
  \ds P_1 ( u, v ) = (1-x)^{-10} (1-2x) \\[2pt]
  \ds P_\infty ( u, v ) = (1-x)^{-4} 
 \end{array}\right .\\
 \left( \frac{7}{10}, \frac{3}{10}, \frac{4}{3} \right) &\left\{ 
 \begin{array}{llll}
  \ds P_0 ( v, u ) = x^{-\frac{19}{3}} (1-x)^\frac{1}{3} \\[2pt]
  \ds P_1 ( v, u ) = x^{-10} (1-2x) \\[2pt]
  \ds P_\infty ( v, u ) = x^{-4} 
 \end{array}\right . &
 \left( \frac{31}{30}, \frac{19}{30}, \frac{4}{3} \right) &\left\{ 
 \begin{array}{llll}
  \ds P_0 ( v, u ) = x^{-\frac{19}{3}} (1-x)^{-\frac{19}{3}} \\[2pt]
  \ds P_1 ( v, u ) = x^{-10} (1-x)^{-10} (1-2x) \\[2pt]
  \ds P_\infty ( v, u ) = x^{-4} (1-x)^{-4} 
 \end{array}\right .\\
\end{array}$$
\par\smallskip\noindent{$\bullet$\quad$I8=\left[ 5, 5, \frac{3}{2} \right] $}
$$ \begin{array}{llll}
  \ds P_0 ( u, v ) := u^{11} v - \frac{11}{432} u^6 v^6 - \frac{1}{186624} u v^{11} \\[10pt]
  \ds P_1 ( u, v ) := u^{30} - \frac{29}{24} u^{25} v^5 - \frac{3335}{62208} u^{20} v^{10} - \frac{3335}{11609505792} u^{10} v^{20} + \frac{29}{835884417024} u^5 v^{25} \\[10pt]
  \ds~~~~~~~~~~~ + \frac{1}{649983722678624} v^{30} \\  
  \ds P_\infty ( u, v ) := u^{20} + \frac{19}{36} u^{15} v^5 + \frac{247}{93312} u^{10} v^{10} - \frac{19}{6718464} u^5 v^{15} + \frac{1}{3482851767} v^{20} \\[10pt]
 \end{array}$$
$$\begin{array}{llll}
 \left( \frac{19}{30}, -\frac{1}{30}, \frac{4}{5} \right) &\left\{ 
 \begin{array}{llll}
  \ds P_0 ( u, v ) = x^\frac{1}{5} ( 1-x )^\frac{1}{5} \\[2pt]
  \ds P_1 ( u, v ) = 1 - 2x \\[2pt]
  \ds P_\infty ( u, v ) = 1 
 \end{array}\right . &
 \left( \frac{1}{6}, \frac{5}{6}, \frac{4}{5} \right) &\left\{ 
 \begin{array}{llll}
  \ds P_0 ( u, v ) = x^\frac{1}{5} ( 1-x )^{ - \frac{11}{5}} \\[2pt]
  \ds P_1 ( u, v ) = (1 - x)^{-6} (1 - 2x) \\[2pt]
  \ds P_\infty ( u, v ) = (1 - x)^{-4} 
 \end{array}\right .\\
 \left( \frac{1}{6}, \frac{5}{6}, \frac{6}{5} \right) &\left\{ 
 \begin{array}{llll}
  \ds P_0 ( v, u ) = x^{-\frac{11}{5}} ( 1-x )^\frac{1}{5} \\[2pt]
  \ds P_1 ( v, u ) = x^{-6} (1 - 2x) \\[2pt]
  \ds P_\infty ( v, u ) = x^{-4} 
 \end{array}\right . &
 \left( \frac{31}{30}, \frac{11}{30}, \frac{6}{5} \right) &\left\{ 
 \begin{array}{llll}
  \ds P_0 ( v, u ) = x^{-\frac{11}{5}} ( 1-x )^{-\frac{11}{5}} \\[2pt]
  \ds P_1 ( v, u ) = x^{-6} (1 - x)^{-6} (1 - 2x) \\[2pt]
  \ds P_\infty ( v, u ) = x^{-4} (1 - x)^{-4} 
 \end{array}\right .\\
\end{array}$$
\par\medskip\noindent
For the following cases, consult \cite{Kob}:
\par\smallskip\noindent{$I9=\left[ 2, 5, \frac{5}{2} \right] $}, {$I10=\left[ 3, 5, \frac{5}{3} \right] $}, {$I11=\left[ \frac{5}{2}, \frac{5}{2}, \frac{5}{2} \right] $}, {$I12= \left[ 3, 5, \frac{3}{2} \right] $}, {$I13=\left[ 5, 5, \frac{5}{4} \right] $}, {$I14=\left[ 2, 3, \frac{5}{2} \right] $}, { $I15= \left[ 3, \frac{5}{2}, \frac{5}{3} \right] $}.
\section{Further study of dihedral cases}
In \S2.3.1, we stated that when $(a,b,c)=(\frac1{2n},-\frac1{2n},\frac12)$, the Kummer solutions $(u,v)$ around $x=0$ satisfy the relation $u^2+\frac1{n^2}v^2=1$. In terms of the Kummer solutions 
$$z^aF(a,a-c+1,a-b+1;z),\quad z^bF(b,b-c+1,b-a+1;z),\quad z=\frac1x,$$
around $x=\infty$,
this relation reads (via the relation connecting the two pairs of Kummer solutions)
$$F\left(\frac{1}{2n},\frac{1}{2n}+\frac{1}{2},\frac{1}{n}+1;{z}\right)F\left(-\frac{1}{2n},-\frac{1}{2n}+\frac{1}{2},-\frac{1}{n}+1;{z}\right)=1,$$
which happens to hold also for non-integral $n$; we write $\frac1n=\alpha$.

Making the same computation for the integer-shifted parameters
$$(a,b,c)=\left(\frac\alpha2-p,-\frac\alpha2-q,\frac12-r\right),\quad p,q,r\in\Z,$$
we get the product
$$z^{-p-q}F\left(\frac{\alpha}{2}{-}p,\frac{\alpha{+}1}{2}{-}p{+}r,\alpha{-}p{+}q{+}1;z\right)F\left(-\frac{\alpha}{2}{-}q,\frac{1{-}\alpha}{2}{-}q{+}r,1{-}\alpha{+}p{-}q;z\right)$$
of the two Kummer solutions around $x=\infty$, which should be a rational function in $z$, since any integer-shift does not affect the monodromy behavior. Replacing $\alpha$ by $\alpha+p+q$, we see that only two integer-parameters are effective. We set
$$s=\frac{p+q}2,\quad t=\frac{p+q}2-r\ \in \frac12\Z \quad (s+t\in\Z),$$
and define the rational function (we write $x$ in place of $z$)
$$f_{s,t}(x):=F\left(\frac{\alpha}{2}-s,\frac{\alpha+1}{2}-t,\alpha+1;x\right)F\left(-\frac{\alpha}{2}-s,\frac{-\alpha+1}{2}-t,-\alpha+1;x\right).$$
In this section, we find an explicit expression of $f_{s,t}(x).$
\begin{lemma}$$
F(a,b,c;x)F(d,e,f;x)=
\sum_{n=0}^{\infty} {}_4 F_3
\left(
\begin{array}{c}
a,\ b,\ 1-n-f,\ -n\\
c,\ 1-n-d,\ 1-n-e
\end{array}
;1\right)\frac{(d,n)(e,n)}{(f,n)(1,n)}  x^n
$$\end{lemma}
Proof:
$$\begin{array}{ll}
\mbox{LHS}
&=\sum_{l=0}^{\infty}\sum_{n+m=l}\frac{(a,n)(b,n)(d,m)(e,m)}{(c,n)(1,n)(f,m)(1,m)}x^l\\[2mm]
&=\sum_{l=0}^{\infty}\sum_{n=0}^{l}\frac{(a,n)(b,n)(d,l-n)(e,l-n)}{(c,n)(1,n)(f,l-n)(1,l-n)}x^l\\[2mm]
&=\sum_{l=0}^{\infty}\frac{(d,l)(e,l)}{(f,l)(1,l)}\sum_{n=0}^{l}\frac{(a,n)(b,n)(f-n+l,n)(l-n+1,n)}{(c,n)(1,n)(d+l-n,n)(e+l-n,n)}x^l\\[2mm]
&=\sum_{l=0}^{\infty}\frac{(d,l)(e,l)}{(f,l)(1,l)}\sum_{n=0}^{l}\frac{(a,n)(b,n)}{(c,n)(1,n)}\textstyle\prod\limits_{k=1}^{n}\frac{(-1+f+l-n+k)(l-n+k)}{(d+l-n+k-1)(e+l-n+k-1)} x^l\\[2mm]
&=\sum_{l=0}^{\infty}\frac{(d,l)(e,l)}{(f,l)(1,l)}\sum_{n=0}^{l}\frac{(a,n)(b,n)}{(c,n)(1,n)}\textstyle\prod\limits_{h=1}^{n}\frac{(f-l+h)(-l+h-1)}{(-d-l+h)(-e-l+h)} x^l\\[2mm]
&=\sum_{l=0}^{\infty}\frac{(d,l)(e,l)}{(f,l)(1,l)}\sum_{n=0}^{\infty}\frac{(a,n)(b,n)(-f,n)(-l,n)}{(c,n)(-d-l+1,n)(-e-l+1,n)(1,n)} x^l=\mbox{ RHS}.
\end{array}$$
Substituting $\,a=\frac{\alpha}{2}-s,\,\,b=\frac{\alpha+1}{2}-t,\,\,c=\alpha+1,\,\,d=-\frac{\alpha}{2}-s,\,\,e=\frac{-\alpha+1}{2}-t,\,\,f=-\alpha+1\,$ we have
\begin{corollary}
\[
f_{s,t}(x)=
\sum_{n=0}^{\infty} {}_4 F_3
\left(
\begin{array}{cccc}
\frac{\alpha}{2}-s,&\frac{\alpha+1}{2}-t,&\alpha-n,\quad-n\\[2mm]
\alpha{+}1,&\frac{\alpha}{2}{+}s{-}n{+}1,&\frac{\alpha+1}{2}{+}t{-}n
\end{array}
;1
\right)
\frac{(-\frac{\alpha}{2}{-}s,n)(-\frac{\alpha{-}1}{2}{-}t,n)}{(-\alpha{+}1,n)(1,n)}x^n.
\]
\end{corollary}
Now we state the result.
\begin{theorem}If $s,t\in\frac12\Z, s+t+1\in\N,\alpha\not\in\Z$, the function 
$$f_{s,t}(x):=F\left(\frac{\alpha}{2}-s,\frac{\alpha+1}{2}-t,\alpha+1;x\right)F\left(-\frac{\alpha}{2}-s,\frac{-\alpha+1}{2}-t,-\alpha+1;x\right)$$
is a polynomial in $x$ of degree $2s$ when $d-t+1\in\N$, and $2t-1$ when $t-s\in\N$\end{theorem}
\subsection{Some lemmas}
\begin{lemma}
For $a\in \R$ and $ n \in \N$, we have
$$(a,n)=0 \quad \Longleftrightarrow \quad 1-a,\,\,a+n \in \N.$$
\end{lemma}
\begin{lemma}If $\alpha \notin \Z$,
\[
i,j\in\N \quad \Longleftrightarrow \quad g(\alpha):= {}_3 F_2
\left(
\begin{array}{cccc}
\frac{\alpha}{2}-i+j+1,\,\,\,\alpha-2i+2,\,\,\,1-i-j\\[2mm]
\frac{\alpha}{2}-i-j+2,\,\quad\alpha-i+j+2
\end{array}
;1
\right)\equiv0.
\]
\end{lemma}
Proof: \par\noindent
[$(\Rightarrow)$]　Possible poles of $g(\alpha)$ are 
$\alpha=0,2,\cdots,2i+2j-4$ and $\alpha=-2j,-2j+1,\cdots,i-j-2$.
\par\smallskip\noindent
[[1]] When $i-j<2$, the poles are simple. So we set
$$g(\alpha)=\sum_{k=1-j}^{[\frac{i-j-1}{2}]} \frac{A_k}{\alpha-2k+1}+\sum_{k=-j}^{[\frac{i-j}{2}]-1} \frac{B_k}{\alpha-2k}+\sum_{k=0}^{i+j-1} \frac{C_k}{\alpha-2k}+D,$$
and prove that $A_k=B_k=C_k=D=0$.
\par\smallskip\noindent
Claim: $A_k=0\,\,\,(1-j \le k \le [\frac{i-j-1}{2}])$
\par\smallskip\noindent
Proof: 
$$\begin{array}{ll}
A_k=&\left[(\alpha-2k+1)\left( \sum_{n=0}^{i-j-2k-1} +\sum_{n=i-j-2k}^{i+j-1} \right)\right]_{\alpha=2k-1}\\[2mm]
&=\left[(\alpha-2k+1) \sum_{n=i-j-2k}^{i+j-1} \frac{(\frac{\alpha}{2}-i+j+1,n)(\alpha-2i+2,n)(-i-j+1,n)}{(\frac{\alpha}{2}-i-j+2,n)(\alpha-i+j+2,n)(1,n)} \right]_{\alpha=2k-1}\\[2mm]
&=\Bigg[\frac{(\frac{\alpha}{2}-i+j+1,i-j-2k)(\alpha-2i+2,i-j-2k)(-i-j+1,i-j-2k)}{(\frac{\alpha}{2}-i-j+2,i-j-2k)(\alpha-i+j+2,i-j-2k-1)(1,i-j-2k)} \\[2mm]
&\qquad \times \left.\sum_{n=0}^{2j+2k-1} \frac{(\frac{\alpha}{2}-2k+1,n)(\alpha-i-j-2k+2,n)(-2j-2k+1,n)}{(\frac{\alpha}{2}-2j-2k+2,n)(\alpha-2k+2,n)(i-j-2k+1,n)}\right]_{\alpha=2k-1}\\[2mm]
&=\frac{(-i+j+k+\frac{1}{2},i-j-2k)(-2i+2k+1,i-j-2k)(-i-j+1,i-j-2k)}{(-i-j+k+\frac{3}{2},i-j-2k)(-i+j+2k+2,i-j-2k-1)(1,i-j-2k)}\\[2mm]
&\quad \times \sum_{n=0}^{2j+2k-1} \frac{(-k+\frac{1}{2},n)(-i-j+1,n)(-2j-2k+1,n)}{(-2j-k+\frac{3}{2},n)(1,n)(i-j-2k+1,n)}\\[2mm]
\end{array}$$
If we write the summand in the last term as  $\sum_{n=0}^{2j+2k-1} a(n)$, then $a(n)+a(2j+2k-n)=0.$ Indeed, we have

$$\begin{array}{ll}
a(2j+2k-n)=&\frac{(-k+\frac{1}{2},2j+2k-n-1)(-i-j+1,2j+2k-n-1)(-2j-2k+1,2j+2k-n-1)}{(-2j-k+\frac{3}{2},2j+2k-n-1)(1,2j+2k-n-1)(i-j-2k+1,2j+2k-n-1)}\\[2mm]
&=\frac{(-k+\frac{1}{2})(-k+\frac{3}{2})\cdots(2j+k-n-\frac{3}{2})}{(-2j-k+\frac{3}{2})(-2j-k+\frac{5}{2})\cdots(k-n-\frac{1}{2})}
\frac{(-i-j+1)(-i-j+2)\cdots(-i+j+2k-n-1)}{(i-j-2k+1)(i-j-2k+2)\cdots(i+j-n-1)}\\[2mm]
&\quad \times\frac{(-2j-2k+1)(-2j-2k+2)\cdots(-n-1)}{1\cdots(2j+2k-n-1)}\\[2mm]
&=\frac{(-k+\frac{1}{2})(-k+\frac{3}{2})\cdots(2j+k-\frac{3}{2})}{(-2j-k+\frac{3}{2})(-2j-k+\frac{5}{2})\cdots(k-\frac{1}{2})}
\frac{(k-n+\frac{1}{2})(k-n+\frac{3}{2})\cdots(k-\frac{1}{2})}{(2j+k-n-\frac{1}{2})(2j+k-n+\frac{1}{2})\cdots(2j+k-\frac{3}{2})}\\[2mm]
&\quad \times\frac{(-i-j+1)(-i-j+2)\cdots(-i+j+2k-1)}{(i-j-2k+1)(i-j-2k+2)\cdots(i+j-1)}\frac{(i+j-n)(i+j-n+1)\cdots(i+j-1)}{(-i+j+2k-n)\cdots(-i+j+2k-1)}\\[2mm]
&\quad \times\frac{(-2j-2k+1)(-2j-2k+2)\cdots(-1)}{1\cdot 2 \cdots(2j+2k+1)}\frac{(2j+2k-n)(2j+2k-n+1)\cdots(2j+2k+1)}{(-n)(-n+1)\cdots(-1)}\\[2mm]
&=(-1)^{2j+2k-1}\frac{(-k+\frac{1}{2})(-k+\frac{3}{2})\cdots(-k+n-\frac{1}{2})}{(-2j-k+\frac{3}{2})(-2j-k+\frac{5}{2})\cdots(-2j-k+n+\frac{1}{2})}\\[2mm]
& \quad \times(-1)^{2j+2k-1}\frac{(-i-j+1)(-i-j+2)\cdots(n-i-j)}{(i-j-2k+1)(i-j-2k+2)\cdots(i-j-2k+n)}\\[2mm]
&\quad\times(-1)^{2j+2k-1}\frac{(-2j-2k+1)(-2j-2k+2)\cdots(-2j-2k+n)}{1 \cdot 2 \cdots n}\\[2mm]
&=-\frac{(-k+\frac{1}{2},n)(-i-j+1,n)(-2j-2k+1,n)}{(-2j-k+\frac{3}{2},n)(i-j-2k+1,n)(1,n)}=-a(n).\\[6mm]
\end{array}$$
\par\smallskip\noindent
Claim: $B_k=0\,\,\,(-j \le k \le [\frac{i-j}{2}]-1)$. 
\par\smallskip\noindent
Proof:
$$\begin{array}{ll}
B_k=&\left[(\alpha-2k)\left( \sum_{n=0}^{i-j-2k-2} +\sum_{n=i-j-2k-1}^{i+j-1} \right)\right]_{\alpha=2k}\\[2mm]
&=\left[(\alpha-2k) \sum_{n=i-j-2k-1}^{i+j-1} \frac{(\frac{\alpha}{2}-i+j+1,n)(\alpha-2i+2,n)(-i-j+1,n)}{(\frac{\alpha}{2}-i-j+2,n)(\alpha-i+j+2,n)(1,n)}\right]_{\alpha=2k}\\[2mm]
&=\Bigg[\frac{(\frac{\alpha}{2}-i+j+1,i-j-2k-1)(\alpha-2i+2,i-j-2k-1)(-i-j+1,i-j-2k-1)}{(\frac{\alpha}{2}-i-j+2,i-j-2k-1)(\alpha-i+j+2,i-j-2k-2)(1,i-j-2k-1)}\\[2mm]
&\qquad \times \left.\sum_{n=0}^{2j+2k} \frac{(\frac{\alpha}{2}-2k,n)(\alpha-i-j-2k+1,n)(-2j-2k,n)}{(\frac{\alpha}{2}-2j-2k+1,n)(\alpha-2k+1,n)(i-j-2k,n)}\right]_{\alpha=2k};\\[2mm]
\end{array}$$
the first term $(-i+j+k+1,i-j-2k-1)$ in the last RHS vanishes.
\par\smallskip\noindent
Claim: $C_k=0\ (0 \le k \le i+j-2)$.
\par\smallskip\noindent Proof:
$$\begin{array}{ll}
C_k
&=\left[(\alpha-2k)\left( \sum_{n=0}^{i+j-k-2} +\sum_{n=i+j-k-1}^{i+j-1} \right)\right]_{\alpha=2k}\\[2mm]
&=\left[(\alpha-2k) \sum_{n=i+j-k-1}^{i+j-1} \frac{(\frac{\alpha}{2}-i+j+1,n)(\alpha-2i+2,n)(-i-j+1,n)}{(\frac{\alpha}{2}-i-j+2,n)(\alpha-i+j+2,n)(1,n)}\right]_{\alpha=2k}\\[2mm]
&=\Bigg[\frac{(\frac{\alpha}{2}-i+j+1,i+j-k-1)(\alpha-2i+2,i+j-k-1)(-i-j+1,i+j-k-1)}{(\frac{\alpha}{2}-i-j+2,i+j-k-2)(\alpha-i+j+2,i+j-k-1)(1,i+j-k-1)}\\[2mm]
&\qquad \times \left.\sum_{n=0}^{k} \frac{(\frac{\alpha}{2}+2j-k,n)(\alpha-i+j-k+1,n)(-k,n)}{(\frac{\alpha}{2}-k+1,n)(\alpha+2j-k+1,n)(i+j-k,n)}\right]_{\alpha=2k};
\end{array}$$
the first term $(-i+j+k+1,i+j-k-1)$ in the last RHS vanishes.
\par\smallskip\noindent
We thus showed that $g(\alpha)$ is a constant $D$, which is equal to
\par\smallskip\noindent$
\left[\sum_{n=0}^{i+j-1} \frac{(\frac{\alpha}{2}-i+j+1,n)(\alpha-2i+2,n)(-i-j+1,n)}{(\frac{\alpha}{2}-i-j+2,n)(\alpha-i+j+2,n)(1,n)}\right]_{\alpha=\infty}
=\sum_{n=0}^{i+j-1}\frac{(-i-j+1,n)}{(1,n)}{\ }{\ }{\ } =(1-1)^{i+j-1}=0.
$
\par\medskip\noindent
[[2]]　When $i-j \ge 2$ the poles at $\alpha = 0,2, \cdots ,i-j-2$ are of order 2, and other poles are of order 1. so we set
$$g(\alpha)=\sum_{k=1-j}^{[\frac{i-j-1}{2}]} \frac{A_k}{\alpha-2k+1}+\sum_{k=-j}^{[\frac{i-j}{2}]-1} \frac{B_k}{\alpha-2k}
+\sum_{k=0}^{[\frac{i-j}{2}]-1} \frac{B'_k}{(\alpha-2k)^2}+\sum_{k=[\frac{i-j}{2}]}^{i+j-1} \frac{C_k}{\alpha-2k}+D,$$
and see that $A_k=B_k=B'_k=C_k=D=0$. We only show
$B'_k=0\,\,\,(0 \le k \le[\frac{i-j}{2}]-1)$; others are quite similar to the previous case.
$$\begin{array}{ll}
B'_k
&=\left[(\alpha-2k)^2\left( \sum_{n=0}^{i+j-k-2} +\sum_{n=i+j-k-1}^{i+j-1} \right)\right]_{\alpha=2k}\\[2mm]
&=\left[(\alpha-2k)^2 \sum_{n=i+j-k-1}^{i+j-1} \frac{(\frac{\alpha}{2}-i+j+1,n)(\alpha-2i+2,n)(-i-j+1,n)}{(\frac{\alpha}{2}-i-j+2,n)(\alpha-i+j+2,n)(1,n)}\right]_{\alpha=2k}\\[2mm]
&=\Bigg[\frac{(\frac{\alpha}{2}-i+j+1,i+j-k-1)(\alpha-2i+2,i+j-k-1)(-i-j+1,i+j-k-1)}{(\frac{\alpha}{2}-i-j+2,i+j-k-2)(\alpha-i+j+2,i-j-2k-2)(\alpha-2k+1,2j+k)(1,i+j-k-1)}\\[2mm]
&\qquad \times \left.\sum_{n=0}^{k} \frac{(\frac{\alpha}{2}+2j-k,n)(\alpha-i+j-k+1,n)(-k,n)}{(\frac{\alpha}{2}-k+1,n)(\alpha+2j-k+1,n)(i+j-k,n)}\right]_{\alpha=2k};
\end{array}$$
the first term $(-i+j+k+1,i+j-k-1)$ in the last RHS vanishes.
\par\smallskip\noindent
[$(\Leftarrow)$]　We assume $g(\alpha) \equiv 0$. The residue $R(-j)$ of $g(\alpha)$ at $\alpha=-j$ is evaluated as follows:

$$\begin{array}{ll}
R(-j)
&=\left[(\alpha+2j)\left( \sum_{n=0}^{i+j-2} +\sum_{n=i+j-1}^{i+j-1} \right)\right]_{\alpha=-2j}\\[2mm]
&=\left[(\alpha+2j) \sum_{n=i+j-1}^{i+j-1} \frac{(\frac{\alpha}{2}-i+j+1,n)(\alpha-2i+2,n)(-i-j+1,n)}{(\frac{\alpha}{2}-i-j+2,n)(\alpha-i+j+2,n)(1,n)}\right]_{\alpha=-2j}\\[2mm]
&=\Bigg[\frac{(\frac{\alpha}{2}-i+j+1,i+j-1)(\alpha-2i+2,i+j-1)(-i-j+1,i+j-1)}{(\frac{\alpha}{2}-i-j+2,i+j-1)(\alpha-i+j+2,i+j-2)(1,i+j-1)}\Bigg]_{\alpha=-2j}\\[2mm]
&= \frac{(-i+1,i+j-1)(-2i-2j+2,i+j-1)(-i-j+1,i+j-1)}{(-i-2j+2,i+j-1)(-i-j+2,i+j-2)(1,i+j-1)}.
\end{array}$$
We have $(-2i-2j+2,i+j-1),(-i-j+1,i+j-1)\neq0$, because if $(-2i-2j+2,i+j-1)=0$ then Lemma 3.4 tells $2i+2j-1,-i-j+1\in \N$, which is a contradiction; if $(-i-j+1,i+j-1)=0$ then Lemma 3.4 tells $i+j+1,0 \in \N$, which is again a contradiction. Thus if $R(-j)=0$ then $(-i+1,i+j-1)=0$; Lemma 3.4 tells $i,j \in \N$.
This proves Lemma 3.5.

\begin{lemma}If $\alpha \notin \Z$,
$(\>\!\textrm{i}\>\!)$\[
s+t+1,\,\,s-t+1\in\N \quad \Longleftrightarrow \quad {}_3 F_2
\left(
\begin{array}{cccc}
\frac{\alpha+1}{2}-t,\,\,\alpha-2s-1,\,\,-2s-1\\[2mm]
\frac{\alpha-1}{2}+t-2s,\,\quad\alpha+1
\end{array}
;1
\right)\equiv0
\]

$(\textrm{ii})$\[
s+t+1,\,\,t-s\in\N \quad \Longleftrightarrow \quad {}_3 F_2
\left(
\begin{array}{cccc}
\frac{\alpha}{2}-s,\,\,\,\,\,\alpha-2t,\,\,\,\,\,-2t\\[2mm]
\frac{\alpha}{2}+s-2t+1,\,\,\alpha+1
\end{array}
;1
\right)\equiv0
\]
\end{lemma}
Proof: Put $i=s+t+1,\,\,j=s-t+1$ in Lemma 3.5, and substitute $\alpha$ by $ \alpha+2t-1$. Then we get (1). The second assertion can be proved similarly.\quad　　　　　　　　　

\subsection{Proof of the theorem}
We suppose $s+t+1 \in \N,\,\, \alpha \notin \Z$.
The function $f_{s,t}(x)$ satisfies the differential equation
$$\begin{array}{ll}
&x^2(1-x)^2w'''+\frac{3}{2}x(1-x)\{2+(2s+2t-3)x\}w'' \\
&+\{(1-\alpha^2)+(\alpha^2-4st+6s+4t-4)x+(2s^2+2t^2+8st-7s-5t+3)x^2\}w' \\
&-\{\alpha^2s+\alpha^2t+2st-s+(4s^2t+4st^2-4st-2s^2+s)x\}w=0,
\end{array}$$
which is the symmetric tensor product of the two hypergeometric equations satisfied by the two hypergeometric series in question. Set
$w=\sum_{n=0}^{\infty} d_nx^n$ and substitute it into the equation above. Then we see that $d_n$ satisfies the difference equation
$$\begin{array}{ll}
&2(n+1)\{\alpha^2-(n+1)^2\}d_{n+1} \\
&-\{4n^3-3(2s+2t-1)n^2-(2\alpha^2-8st+6s+2t-1)n+2(\alpha^2s+\alpha^2t+2st-s)\}d_{n} \\
&+(n-2s-1)(n-2t)(2n-2s-2t-1)d_{n-1}=0.
  \end{array}$$
If we put $n=2s+1$ and $n=2t$, we have
$$d_{2s+2}=\frac{s-t+1}{2s+2}d_{2s+1},\quad \mbox{and}\quad
d_{2t+1}=\frac{s-t}{2t+1}d_{2t},$$
respectively.
\par\noindent
When $s-t+1 \in \N$, since we have $\alpha/2-s=\alpha/2+s-n+1$, in the expression in Corollary 3.2, ${}_4F_3$ reduces to ${}_3F_2$:

\[
d_{2s+1}=
{}_3 F_2
\left(
\begin{array}{cccc}
\frac{\alpha+1}{2}-t,&\alpha-2s-1,\ -2s-1\\[2mm]
\alpha+1,&\frac{\alpha-1}{2}+t-2s
\end{array}
;1
\right)
\frac{(-\frac{\alpha}{2}-s,2s+1)(-\frac{\alpha-1}{2}-t,2s+1)}{(-\alpha+1,2s+1)(1,2s+1)}. 
\]
Thus by Lemma 3.5 we have 
$d_{2s+1}=0$, and so the above relation leads to $d_{2s+2}=0$. Now the difference equation above asserts
$$d_n=0\qquad(n \geq 2s+1).$$
On the other hand, we have $d_{2s}\neq0$, otherwise all the coefficients would vanish. We therefore showed that $\deg f_{s,t}(x)=2s$.
\par\smallskip\noindent
When $t-s \in \N$, the assertion will be similarly proved.

\section{When the projective monodromy group is a Fuchsian group}
We study the affine Schwarz map when the projective monodromy group $\bar G$ is (conjugate to) a Fuchsian group, and the inverse of the Schwarz map is single-valued. The Fuchsian group is of genus zero, and $sch^{-1}$ is a $\bar G$-modular function often called the Hauptmodul.
\par\medskip\noindent
Typical cases: $$\begin{array}{ccc}
 (a,b,c)&\mbox{projective monodromy group}\\[2mm]
(1/12,5/12,1)&PSL_2(\Z)\\[2mm]
(1/2,1/2,1)&\Gamma(2)\quad (\mbox{congruence subgroup})\\[2mm]
(13/84,1/42,1/2)&\mbox{with smallest co-volume}\end{array}$$


\begin{remark}
When $(a,b,c)=(1/2,1/2,1)$, since $\bar G\cong\Gamma(2)$ is a free group, the map $asch$ is defined on the universal covering of $X$, which is isomorphic to the unit disc $\D=\{x\in\C\mid |x|<1\}$. Thus the map $asch$ gives a closed embedding of $\D$ into $\C^2$. Existence of such an embedding was firstly proved by Nishino \cite{Nno}.
\end{remark}
We can assume that the image $sch(X)$ is open dense in the upper half plane $\HH:=\{\tau\in Z\mid \Im\tau>0\}$. Then the curve $C$ lies in the cone $\{(u,v)\in W\mid u/v\in \HH\}$, which we parametrize as 
$$u/v=\tau\in\HH,\quad v\in\C^\times.$$
The curve $C$ is the graph of the function 
$$v=F(a,b,c;sch^{-1}(\tau)),\quad \tau\in\HH,$$
which is a modular form.
When $(a,b,c)=(1/12,5/12,1)$, the inverse of $\tau=sch(x)=u(x)/v(x)$ is given by the elliptic modular function $j(\tau)$, and
$$F\left(\frac1{12},\frac5{12},1;\frac{12^3}{j(\tau)}\right)$$
is a modular form of which fourth power is the Eisenstein series $E_4(\tau)$. 
Such functions are systematically studied by Koike \cite{Koi}.
It might be rare that a modular form is considered as a {\it map}.
\par\smallskip\noindent
{\bf Aknowledgement:} The authors are grateful to T. Sasaki for his encouragement and useful comments.

\bigskip
\begin{flushleft}
%
Department of Mathematics Kyushu University, 
Fukuoka 810-8560 Japan
\end{flushleft}
\end{document}